\documentclass{amsart}
\usepackage{ifthen}
\usepackage{hyperref}

\usepackage{amsmath,amssymb,amsthm,verbatim,paralist,graphicx}

\usepackage[active]{srcltx}
\usepackage[T1]{fontenc}
\usepackage[latin1]{inputenc}
\usepackage{yfonts}
\usepackage{color}
\DeclareMathOperator{\rk}{r}
\DeclareMathOperator{\ind}{ind}
\DeclareMathOperator{\cl}{cl}
\DeclareMathOperator{\supp}{supp}
\DeclareMathOperator{\sep}{sep}
\DeclareMathOperator{\Dir}{Dir}

\numberwithin{equation}{section}
\newtheorem{theorem}{Theorem}[section]
\newtheorem{prop}{Proposition}[section]
\newtheorem{lemma}{Lemma}[section]
\newtheorem{corollary}{Corollary}[section]
\newtheorem{definition}{Definition}[section]

\newtheorem{remark}[lemma]{Remark}

\newtheoremstyle{TheoremNum}
        {\topsep}{\topsep}              
        {\itshape}                      
        {}                              
        {\bfseries}                     
        {.}                             
        { }                             
        {\thmname{#1}\thmnote{ \bfseries #3}}
    \theoremstyle{TheoremNum}
    \newtheorem{thmn}{Theorem}

\def\Oscr{\mathcal{O}}

\begin{document}
\title{Vertex-shellings of Euclidean Oriented Matroids}  \author{Winfried Hochst\"attler} \author{Michael Wilhelmi}
\address{FernUniversit\"at in Hagen}
\email{Winfried.Hochstaettler@FernUniversitaet-Hagen.de}
\email{mail@michaelwilhelmi.de}

\begin{abstract}
  A Euclidean oriented matroid program yields a partial ordering of the cocircuits of its cocircuit graph.
  We show that every linear extension of that ordering yields 
  a topological sweep and
  induces a recursive atom-ordering (a shelling of the cocircuits) of
  the tope cell of the feasible region. We extend that sweep and obtain
  also a vertex-shelling of the whole oriented matroid and finally
  describe some connections to the notion of stackable zontopal tilings and to a counterexample of a conjecture of A. Mandel.
\end{abstract}

\dedicatory{}

\maketitle


\section{Introduction}
{\em Shellability} is a concept from topological combinatorics. Most prominently, it was used by Bruggesser and
Mani~\cite{8} to complete a proof of Schl\"afli for the
Euler-Poincar\'e relation and by McMullen~\cite{McMullen} to prove the
Upper Bound Conjecture.  A {\em shelling} or a {\em recursive coatom
  ordering} of a graded poset is an ordering of coatoms that can be
made compatible with a recursive coatom ordering of the ideal induced
by each coatom.  To show that the boundary complex of a convex
polytope is shellable Bruggesser and Mani in \cite{8} used the concept
of a {\em line shelling}.  This easily generalizes to oriented
matroids as follows. If $P$ is a tope and $\ell$ is a shortest path from $P$ to -$P$
in the tope graph, then recording the facets of $P$, in the order their
defining hyperplanes are met by $\ell$, yields a shelling order of the tope.

While the dual of a face lattice of a polyhedron is again
polyhedral, this is no longer true for oriented matroids (\see e.g.\
\cite{2}, 9.3.10). It is an open problem, whether duals of face
lattices of oriented matroid polytopes admit a shelling, i.e.\ if
topes admit a recursive atom ordering (\cite{2}, Exercise 4.18 (b)).
The ``polar procedure'' of tracing a line yields a hyperplane sweep,
and in polytopes this yields an ordering of the vertices
corresponding to a recursive atom ordering.

In non-realizable oriented matroids such a topological sweep is not
always possible (\cite{2}, 10.4.6), namely, in the oriented matroid
generalizations of linear programs we may have cycles where each edge
``properly improves'' the objective function. Any hyperplane sweep
must entangle itself in such a cycle.  Oriented matroid programs, where
such "improving cycles" do not occur, are called {\em Euclidean}.
It is mentioned in \cite{2}, after 10.4.6, and made explicit in
\cite{1} that topological sweeping is possible in Euclidean oriented
matroid programs. 

It seems immediate that such a topological sweep yields a recursive
atom ordering of the feasible region of the oriented matroid program,
but the recursive part turns out to be suprisingly difficult. While it is
clear that a minor of a Euclidean oriented matroid program still is
Euclidean, see \cite{5}, 9.III.17 or \cite{2}, 10.5.6, we may lose
the boundedness, which is necessary for our shelling, in some of the
minors to be considered.  We fix this by using lexicographic extensions 
to bound the minors, and in our companion paper~\cite{Paper1} we show that these extensions preserve Euclideaness. 

In Section~\ref{sec:two} of this paper we present some preliminaries and
required definitions, 
Section~\ref{sec:five} reviews the
topological sweep and shows the following theorem:
\begin{theorem}\label{theo:topSweepExists}
  Let $(\Oscr,g,f)$ be a Euclidean oriented matroid program. Then $G_f$ defines a partial order (poset) of the cocircuits of $(\Oscr \setminus f)$.
  Let $f$ be in general position.
  Every linear extension of that poset defines a non-degenerate topological sweep and vice versa.
\end{theorem}

The difficult part is to show that a linear extension yields a sweep. We show it by single-element extensions 
of the matroid maintaining the order of the linear extension.    
In Section \ref{sec:four} we show the vertex-shelling of the tope lattice and the whole oriented matroid. 

\begin{theorem} \label{theo:FourthMainTheoem}
  Let $\mathcal{O}$ be a Euclidean oriented matroid with $f \in
  \mathcal{O}$ lying in general position.  Then $\mathcal{O} \setminus f$ and each tope cell of
  $\mathcal{O} \setminus f$ are vertex-shellable.
\end{theorem}

We end with our conclusion, relating our work to recent research, e.g.\
the stackability of zonotopal tilings, see \cite{Athanasiadis}  and presenting some
open questions and the relation to the work of Mandel.

\section{Preliminaries}\label{sec:two}

This paper  uses the results of our companion paper~\cite{Paper1}. 
We cite here the required results. omitting the proofs.
For the basic concepts/notations of oriented matroid programming, we
refer to \cite{2}, Chapter 10. Unlike in \cite{2}, we call an oriented
matroid program $(\mathcal{O},g,f)$ {\em bounded} if $X_g = +$ for
every cocircuit $X$ of the feasible region. 
Using ``complementary
slackness'' and the Main Theorem of Oriented Matroid Programming
\cite{2}, 10.1.12/13 it is immediate that a feasible, bounded oriented
matroid program has at least one maximal and one minimal solution.  
We
use the notion of Euclideaness of an oriented matroid program in
\cite{2}, 10.5.2. An oriented matroid program
  $(\mathcal{O},g,f)$ is Euclidean if and only if its cocircuit graph
  $G_f$ contains no directed cycles.  Additionally we call an
oriented matroid $\mathcal{O}$ {\em Euclidean} if all programs
$(\mathcal{O},g,f)$ with $f \neq g$ of the groundset ($f$ not a
coloop, $g$ not a loop) are Euclidean and we call it {\em Mandel} if
there exists a single element extension $\mathcal{O}'
= (\mathcal{O} \cup f)$ of $\mathcal{O}$ in general position such that
$(\mathcal{O}',g,f)$ is Euclidean for all elements $g$ of the
groundset that are not loops.  We call an oriented matroid
$\mathcal{O}$ {\em strongly Euclidean} if it is of rank $1$ or possesses
an element $g$ (not a loop) such that $\mathcal{O} / g$ is strongly
Euclidean and the program $(\mathcal{O}'',g,f)$ is Euclidean for all
extensions $\mathcal{O}''=\mathcal{O} \cup f$.

\subsection{Directing Edges in the Graph $G_f$ and Beyond}

Recall, that if $X$ and $Y$ are covectors in an oriented matroid
$\mathcal{O}$ we have $z(X \circ Y) = z(Y \circ X) = z((-X) \circ Y)$
for their zero sets. Moreover, if $X,Y$ are cocircuits with $X \neq
\pm Y$ then $\supp(Y) \setminus \supp(X) \ne \emptyset \ne \supp(X)
\setminus \supp(Y)$. A covector $F$ is an {\em edge} in $\mathcal{O}$
if $z(F)$ is a coline in the underlying matroid
$\mathcal{M}(\mathcal{O})$.  We say a cocircuit $X$ is {\em lying in the edge $F$ }or {\em on the line spanned by} (or {\em through}) $F$ if $z(X) \supset z(F)$. We call two covectors a {\em modular
  pair} or {\em comodular} in $\mathcal{O}$ if their zero-sets form a
modular pair in $\mathcal{M}(\mathcal{O})$.  Hence two cocircuits $X$
and $Y$ are comodular if and only if $X \circ Y$ is an edge in~$\mathcal{O}$.  Also we will call $(X,Y)$ an {\em
  edge} in $\mathcal{O}$ if $X$ and $Y$ are conformal comodular cocircuits.

The following propositions are immediate.

\begin{prop}[\cite{Paper1}]\label{prop:CocircuitsLyingOnEdges1}
Let $F$ be an edge in an oriented matroid $\mathcal{O}$ and
$X \neq \pm Y$ cocircuits and  $z(F) \subseteq z(X)$ as well as  $z(F) \subseteq z(Y)$. Then
$z(X \circ Y) = z(F).$
\end{prop}

\begin{prop}[\cite{Paper1}]\label{prop:unqieCocircuitOnEndge}
  Let $F$ be an edge of an oriented matroid $\mathcal{O}$ and $e \in
  \supp(F)$.  Then there is, up to sign reversal, a unique cocircuit
  $Z$ with $e \not \in \supp(Z) \subseteq \supp(F)$.
\end{prop}

\begin{prop}[Uniqueness Of Cocircuit Elimination, see \cite{Paper1}]\label{prop:UniqueCoEl}
  Let $X,Y$ be comodular cocircuits of an oriented matroid and $e \in
  \sep(X,Y)$. Then cocircuit elimination of $e$ between $X$ and $Y$
  yields a unique cocircuit $Z$ such that $\supp(Z) \subseteq \supp(X
  \circ Y) \setminus e$. Furthermore, we have $Z_f=X_f\circ Y _f$ for
  all $f \notin \sep(X,Y)$.
\end{prop}

We introduce a direction for modular pairs of
cocircuits.

\begin{definition}[\cite{Paper1}]\label{def:Dirfunction}
Let $(\mathcal{O},g,f)$ be an oriented matroid program.
Let 
\[ \mathcal{P} = \{ (X,Y) \mid X \text{ and } Y \text{ are comodular cocircuits in } \mathcal{O} \text{ with } X_g = Y_g \neq 0 \}.\]

We say $(X,Y) \in \mathcal{P}$ is directed from $X$ to $Y$ and write
$X \to_{g,f} Y$ \\ (or $X \to_{(\mathcal{O},g,f)} Y$), if cocircuit
elimination of $g$ between $-X$ and $Y$ yields $Z$ satisfying $Z_f=+$.
We write $X \leftarrow_{g,f} Y$ if $Z_f=-$ and $X
\leftrightarrow_{g,f} Y$ if $Z_f=0$ (then we say $(X,Y)$ is
\em{undirected}).  We may omit the subscript, in case it is obvious.

\end{definition} 

The following again is immediate (we use the order $- \prec 0 \prec +$ in $\{-,0,+\}$):

\begin{prop}[\cite{Paper1}] \label{prop:basicDirectionProps} Let $(X,Y) \in \mathcal{P}$. Then it holds
  \begin{enumerate} 
\item $X \to Y \iff -Y \to -X$.
   \item  $X \to_{(\mathcal{O},g,f)} Y \Leftrightarrow X \to_{({}_{-g}\mathcal{O},g,f)} Y$ and $X \to_{(\mathcal{O},g,f)} Y \Leftrightarrow Y \to_{({}_{-f}\mathcal{O},g,f)} X$ where ${}_{-g}\mathcal{O}$ resp.  ${}_{-f}\mathcal{O}$ are derived from $\mathcal{O}$ reorienting $g$ resp. $f$.
    \item If $X_f \prec Y_f$, then $X \to
    Y$. If $X_f = Y_f = 0$, then $X \leftrightarrow Y$.

  \item Let $e \in \sep(X,Y)$, $X \to Y$
    and cocircuit elimination of $e$ between $X$ and $Y$ yield $W$.
    Then $(X,W),(W,Y) \in \mathcal{P}$ and 
    $X \to W$ holds as well as $W \to Y$.
   \item If $(X,Y)$ is undirected, then also $(Y,X), (-X,-Y)$ and $(-Y,-X)$ are undirected in every reorientation of $\mathcal{O}$.
  \end{enumerate}
\end{prop}

We add a small Proposition.
\begin{prop}[\cite{Paper1}] \label{prop:TransitivityOfDir}
Let $X_g = Y_g \neq 0$ and $X_f = Y_f \neq 0$ and $h \in E \setminus \{f,g\}$. If $X \leftrightarrow_{(\mathcal{O},g,f)} Y$, then  $X \to_{(\mathcal{O},g,h)} Y$ implies $X \to_{(\mathcal{O},f,h)} Y$.
\end{prop} 
\begin{proof}
Cocircuit elimination of $g$ between $-X$ and $Y$ yields $Z$ with $Z_f = 0$. Hence, elimination of $f$ between $-X$ and $Y$ yields again $Z$.
\end{proof}

If $X$ and $Y$ are a conformal comodular pair, the former
definition is compatible with the orientation of the cocircuit graph
of an oriented matroid program.

\begin{definition} [see \cite{2}, 10.1.16]\label{def:Graph}
Let $(\mathcal{O},g,f)$ be an oriented matroid program. 
Let $G_f$ be the graph whose vertices are the cocircuits of $(\mathcal{O},g)$ and whose edges are the edges $(Y^1,Y^2)$ between two comodular, conformal cocircuits.
The edges in $G_f$ are (un-)directed like in $(\mathcal{O},g,f)$.
We call the graph $(G \setminus f)_f$ the directed subgraph of $G_f$
induced by the cocircuits of $(\mathcal{O} \setminus f,g)$.  If two
comodular cocircuits $Y^1,Y^2$ with $Y^1 \circ Y^2$ are conformal in
$(\mathcal{O} \setminus f,g)$ but not in $G_f$, then we have $Y^1_f \succ
Y^2_f$ (or $Y^1_f \prec Y^2_f$) in $\mathcal{O}$ and we direct the edge 
$Y^1 \rightarrow Y^2$ (or $Y^1 \leftarrow Y^2$) in $(G \setminus f)_f$.
\end{definition}

We define paths in $G_f$. We give a reformulation of Definition 10.5.4 in \cite{2}.

\begin{definition}[\cite{Paper1}]
A {\em path} in $G_f$ from $P^1$ to $P^n$ is a sequence
\[ P = [P^1, \hdots, P^n] \]
of vertices of $G_f$ such that $(P^i,P^{i+1})$ is an edge in $G_f$ for all $1 \leq i < n$.
A path $P$ is {\em undirected} if $P^i \leftrightarrow P^{i+1}$ for
all $1 \leq i < n$ and {\em directed} if is not undirected and either
$P^i \to P^{i+1}$ or  $P^i \leftrightarrow P^{i+1}$ for all $1
\leq i < n$, in which case we say that the path {\em is
  directed from $P^1$ to $P^n$}.

A path is {\em closed} if $P^1 = P^n$. A closed path is a {\em cycle}
if $P^1, \hdots, P^{n-1}$ are pairwise distinct.
\end{definition}

We collect some more well-known facts.

\begin{prop}[\cite{Paper1}]  \label{prop:PathGraphWellKnownFacts}
  \begin{enumerate}
  \item \label{prop:ClosedPathContainsCycle}
Let $P = [P^1, \hdots, P^n]$ be a closed path in a graph and and let $X,Y$ be subsequent vertices of $P$. 
Then there exists always a subpath $Q$ of $P$ containing $X$ and $Y$ that is a cycle.
\item \label{pathsOnEdges} If $X,Y$ are comodular cocircuits in $G_f$
  then there is a path $X^0, \hdots, X^n$,  $X=X_0, Y=X_n$ in $G_f$ with
  distinct cocircuits such that $\supp(X_i \circ X_j) = \supp(X \circ
  Y)$ and all pairs $(X_i, X_j)$ are directed as $(X,Y)$ for all $i < j$.
\item \label{prop:DirectedPaths}
Let $P$ be a directed path $P = (X,X^1, \hdots, X^n, Y)$ between two cocircuits $X$ and $Y$ in $G_f$.
If $X_f \prec Y_f$ the path 
is directed
from $X$ to $Y$.
\item \label{prop:PathGraphWellKnownFacts_directedCycle} Let $(\mathcal{O},g,f)$ be an
  oriented matroid program that contains a directed cycle in the graph
  $G_f$. Then all cocircuits in the cycle have the same $f$-value, in particular $X_f \ne 0$ holds everywhere.
  \end{enumerate}
\end{prop}

%
%
%
%
%
%
%

In Euclidean oriented matroid programs exchanging the roles of infinity and the target function preserves Euclideaness.

\begin{theorem}[\cite{Paper1}] \label{theorem:EuclideanessStays}
If $(\mathcal{O},g,f)$ is a Euclidean oriented matroid program, so is $(\mathcal{O},f,g)$.
\end{theorem}

Clearly, reorienting elements from $E \setminus \{g,f\}$ does not
change the direction of any edge, and hence:

\begin{lemma}[see \cite{Paper1}, Lemma 2.3]
Let $(\mathcal{O},g,f)$ be a Euclidean oriented matroid program with groundset $E$ and let
$S \subseteq E \setminus \{f,g\}$ and let  ${}_{-S}\mathcal{O}$ be derived from $\mathcal{O}$ by reorientation of the elements of $S$.
Then $({}_{-S}\mathcal{O},g,f)$ is a Euclidean oriented matroid program, too.
\end{lemma}
 
\subsection{Old and new cocircuits in single-element extensions} 

For the notion of single-element extensions of oriented matroids and its corresponding localizations, we refer to \cite{2}, Chapter 7.1.
We call an extension non-trivial if it is not an extension with a coloop. We mention one fact explicitly.

\begin{theorem}[see \cite{2}, 7.1.4]\label{theorem:ext}
Let $\mathcal{O}$ be an oriented matroid of rank $\rk$ and $\mathcal{C}$ the set of its cocircuits,
let $\mathcal{O}' = \mathcal{O} \cup p$ be a non-trivial 
single-element extension of $\mathcal{O}$
and let $\sigma$ be the corresponding localization. Let 
\[  S =  \{Z = (Y,\sigma(Y)) \colon Y \in \mathcal{C} \}.  \] 
Then the cocircuits $ \mathcal{C}'$ of $\mathcal{O}'$ are given by
\begin{align*} \mathcal{C}' = S \cup \{ Z' = (Y^1 \circ Y^2,0) \colon & Y^1, Y^2 \in \mathcal{C}, \sigma(Y^1) = -\sigma(Y^2) \neq 0,  \\
 & \sep(Y^1,Y^2) = \emptyset, \rk(z(Y^1 \circ Y^2)) = \rk -2 \} 
 \end{align*}
 We say that the cocircuit $Z \in S$ is {\em derived from} the
 cocircuit $Y$ or {\em old} and the cocircuit $Z' \in \mathcal{C}'
 \setminus S$ is {\em derived from the edge} $Y^1 \circ Y^2$ or {\em
   from the cocircuits} $Y^1,Y^2$ or {\em new}. We call an edge of two
 old cocircuits old as well. We call an edge new if at least one of its cocircuits is new.
 \end{theorem}
 Let us mention some well-known relations between edges and vertices
 in the oriented matroid and its extension.
 \begin{prop}[see \cite{Paper1}, Proposition 2.8] \label{prop:uniqueEdges2} Let $\mathcal{O}$ be an
   oriented matroid and let $\mathcal{O}' = \mathcal{O} \cup p$ be a
 non-trivial single-element extension of $\mathcal{O}$.
  \begin{enumerate}
  \item \label{prop:uniqueEdges} 
    If $Z = (X \circ Y,0)$ is a cocircuit in $\mathcal{O}'$
    derived from the edge $X \circ Y$ where $X,Y$ are conformal
    comodular cocircuits in $\mathcal{O}$, then the pair $\{X,Y\}$ is
    unique.
  \item \label{prop:cocircuitsInSingleElementExtensions} If $Y$ is a
    cocircuit of $\mathcal{O}'$ with $Y_p \neq 0$, then $Y \setminus
    p$ is a cocircuit in $\mathcal{O}$, $Y$ is old.
\item If $p$ is in general position, then every cocircuit $X$ of $\mathcal{O}'$ with $X_p = 0$ is new.
  \end{enumerate}
\end{prop}

\subsection{Lexicographic extensions}

We start with a technical definition.

\begin{definition}[\cite{Paper1}]
  Let $I = [e_1, \hdots, e_k ]$ be an ordered subset of the groundset
  of an oriented matroid $\mathcal{O}$ and let $Y$ be a cocircuit in
  $\mathcal{O}$. We call the index $i$ such that $Y_{e_i} \neq 0$ and
  $Y_{e_j} = 0$ for all $j < i$ the {\em index of the cocircuit $Y$
    wrt.\ I}. If $Y_{e_i} = 0$ for all $1 \leq i \leq
  k$ let the index of the cocircuit wrt.\ I be $k+1$.
\end{definition}

We define a lexicographic extension as follows.

\begin{definition}[see \cite{2}, 7.2.4] \label{prop:lexExtBjoerner}
Let $\mathcal{O}$ be an oriented matroid with  $\mathcal{C}$ being its set of cocircuits, $I = [e_1, \hdots, e_k]$ an ordered subset of the groundset, and $\alpha = [\alpha_1, \hdots, \alpha_k] \in \{+,-\}^k.$
Then the lexicographic extension $\mathcal{O}[I^{\alpha}] = \mathcal{O}[e^{\alpha_1}_1, \hdots, e^{\alpha_k}_k]$ of $\mathcal{O}$ is given by the localization
$\sigma \colon \mathcal{C} \rightarrow \{+,-,0\}$ with
\[ \sigma(Y) = \begin{cases}\alpha_iY_{e_i} & \text{ if for the index } i \text{ of }Y \text{ wrt.\ } I \text{ holds } i \leq k , \\ 0 &  otherwise. \end{cases}\]
We call the lexicographic extension \em{positive} and write $\mathcal{O}[I^+]$ iff $\alpha_i = +  $ for all $i \leq k$.
\end{definition}

It is immediate that a lexicographic extension can always be written
as a reorientation of the original, a positive lexicographic extension
and then reversing the reorientation in the extension.  Until the end
of this chapter, let $\mathcal{O}$ always be an oriented matroid of
rank $\rk \geq k$, let $I = [e_1, \hdots, e_k]$ be an ordered subset of
$\mathcal{O}$, let $\alpha = [\alpha_1, \hdots, \alpha_k] \in
\{+,-\}^k$  and let $\mathcal{O}' = \mathcal{O}[I^{\alpha}] =
(\mathcal{O} \cup p)$ be the lexicographic extension of $\mathcal{O}$.

\begin{prop}[see \cite{Paper1}, Proposition 2.9]\label{prop:IndexPUnequalZero}
The index $\ind_Y$ wrt.\ $I$ for a cocircuit $Y$ in $\mathcal{O}[I^{\alpha}]$ with $Y_p \neq 0$  is always at most $k$ and 
$Y_{e_{\ind_Y}} = Y_p$.  
\end{prop}

Recall that we call an extension {\em principal} if we add a point to a flat in general 
position in the matroid.  

\begin{prop}[see \cite{Paper1}, Proposition 2.10]\label{prop:lexExtensionIsFree}
  The lexicographic extension $\mathcal{O}[I^{\alpha}]$ is a principal
  extension of $\mathcal{O}$ adding a point to the flat $\cl(e_1,
  \hdots, e_k)$.
\end{prop}


%
%

  We can always extend an oriented matroid program to a bounded one.

  \begin{lemma}[see \cite{Paper1}, Lemma 2.6]\label{lem:LexExtensionIsBounded} Let
    $(\mathcal{O},g,f)$ be an oriented matroid program.  Let
    $\mathcal{O}' = \mathcal{O} \cup p = \mathcal{O}[I]^+$ be a
    positive lexicographic extension of $\mathcal{O}$ in general
    position where $I$ is a base of $\mathcal{O} \setminus f$. Then
    $(\mathcal{O}',p,f)$ is bounded. It is feasible if and only if
    $(\mathcal{O},g,f)$ is.  \end{lemma} 

  \subsection{Directions of edges of old cocircuits in
    single-element extensions}\label{subs:dirOfEdgesOfOld} We start
  examining the directions of old edges.  

    \begin{prop}[see \cite{Paper1}, Proposition 3.1]\label{prop:EdgeStaysEdge} Let $(\mathcal{O},g,f)$ be
      an oriented matroid program and let $\mathcal{O}' = \mathcal{O}
      \cup p$ be a non-trivial single-element extension of
      $\mathcal{O}$.  Let $X,Y$ be two old conformal, comodular
      cocircuits in $(\mathcal{O}',g)$.  Then $X \setminus p,Y
      \setminus p$ are conformal cocircuits in $(\mathcal{O},g)$.
      They are also comodular in $(\mathcal{O},g)$ if $(X \circ Y)_p
      \ne 0$ or if $\mathcal{O} \cup p$ is a principal extension (or
      if both cases hold).  \end{prop} 
      
  \begin{remark}[see \cite{Paper1}, Remark 3.1]  If $(X \circ Y)_p = 0$ and
      $\mathcal{O} \cup p$ is not a principal extension it is possible
      that $X \setminus p$ and $Y \setminus p$ are not comodular in
      $\mathcal{O}$ e.g. if $\mathcal{O}$ has rank 3 and $z(X
      \setminus p)$ and $z(Y \setminus p)$ are disjoint lines in the
      underlying matroid intersecting at the point $p$ in the
      extension.  Because the lexicographic extension is principal
      (see Proposition \ref{prop:lexExtensionIsFree}) the two
      cocircuits are always comodular in $\mathcal{O}$ in that case.
     \end{remark} 

We introduce a new function to    abbreviate our notions.  

\begin{definition}[Dir-function, see \cite{Paper1}, Definition 3.1] Let
    $X,Y$ be a comodular pair of cocircuits in $(\mathcal{O},g,f)$
    with $X_g = Y_g \neq 0$. Then we define
    \[\Dir_{(\mathcal{O},g,f)}(X,Y) = \left\{ \begin{array}[h]{lcr}
        + & \text{ iff \quad} X \rightarrow Y,\\
        - & \text{ iff \quad} X \leftarrow Y,\\
        0 & \text{ iff \quad} X \leftrightarrow Y.  \end{array}
    \right\} \text{ in } (\mathcal{O},g,f).\] If $X,Y$ are comodular
    cocircuits with $X_g \neq 0$ and $Y_g = 0$ (or vice versa) then we
    define $\Dir(X,Y) = Y_f$ (or $\Dir(X,Y) = -X_f$ resp.).
  \end{definition} 
 
Finally, we remark that the directions of old edges are
  preserved in extensions.  

\begin{lemma}[direction preservation of
    old edges, see \cite{Paper1}, Lemma 3.2]\label{lem:DirPresOfOutsideEdges} Let
      $(\mathcal{O},g,f)$ be an oriented matroid program and let
      $\mathcal{O}' = \mathcal{O} \cup p$ be a non-trivial
      single-element extension of $\mathcal{O}$.  Let $X,Y$ be two
    old cocircuits in $(\mathcal{O}',g)$ with $X \circ Y$ being an
    edge in $\mathcal{O}'$.  If $((X \circ Y) \setminus p)$ is also an
    edge in $\mathcal{O}$ then cocircuit elimination of $g$ between
    $-X$ and $Y$ yields an old cocircuit and it holds \[
    \Dir_{(\mathcal{O}',g,f)}(X,Y) =
      \Dir_{(\mathcal{O},g,f)}(X \setminus p, Y \setminus p). \]
  \end{lemma} 
  
 \subsection{The directions of the cocircuit graph of the
    lexicographic extension $(\mathcal{O} \cup
    p,p,f)$}\label{section:DirGraphpf} We can describe completely the
  directions of the graph $G_f$ of $(\mathcal{O} \cup p,p,f)$ if we
  know the directions of some edges in the graphs $G_f$ of
  $\mathcal{O}$.  Until the end of this chapter, let $\mathcal{O}$ be
  an oriented matroid of rank $\rk$ with groundset $E$ and let $f \in
  E$.  Let ${I} = [e_1, \hdots, e_k]$ with $k \leq \rk$ be an ordered
  set of independent elements of $E$, where possibly $f \in I$.  Let
  $\mathcal{O}' = \mathcal{O} \cup p$ be the lexicographic extension
  $\mathcal{O}[e^+_1, \hdots, e^+_k]$.  
We need two technical lemmas.
   
\begin{lemma}[see \cite{Paper1}, Lemma 3.3]\label{lem:case3lexext}
  Let $X,Y$ be two conformal
    cocircuits with $X_p = Y_p \neq 0$ and $X \circ Y$ being an edge
    in $\mathcal{O}'$.  
    Let $\ind_X, \ind_Y$ be
    the indices of $X, Y$ wrt.\ $I$. 
  Let $1 \leq i = \ind_X = \ind_Y \leq k$ and
  $f \neq e_i$ and thus $X_{e_i} = Y_{e_i}\neq 0$.  If $X_f \neq Y_f$
  or $X_f = Y_f = 0$ or  $Dir_{e_i,f}(X,Y) \neq 0$ then 
\[\Dir_{e_i,f}(X,Y) = \Dir_{p,f}(X,Y) \text{ and } \Dir_{f,e_i}(X,Y) = \Dir_{f,p}(X,Y). \]
\end{lemma}

\begin{lemma}[see \cite{Paper1}, Lemma 3.5]\label{lem:case6lexExt}
  Let $X,Y$ be two conformal
    cocircuits with $X_p = Y_p \neq 0$ and $X \circ Y$ being an edge
    in $\mathcal{O}'$.  
    Let $\ind_X, \ind_Y$ be
    the indices of $X, Y$ wrt.\ $I$. 
Let $1 \leq i = \ind_X < \ind_Y = j \leq k$ and let $Y_f \neq 0$. Then 
\[\Dir_{p,f}(X,Y) = \Dir_{e_i,f}(X,Y) = Y_f. \]
\end{lemma}

%
%
We cite the main theorem of our companion paper [\cite{Paper1}].

\begin{theorem}[see \cite{Paper1}, Theorem 1.2]\label{Theo:SecondMaintheorem2}
A lexicographic extension of a Euclidean oriented matroid
is  Euclidean.
\end{theorem}

Here we need it in a more special form.

\begin{corollary}[see \cite{Paper1}, Corollary 3.1]\label{cor:lexExtStaysEucl2}
Let $\mathcal{O}$ be an oriented matroid of rank $\rk$ with groundset $E$ and let ${f \in E}$ being in general position. 
Let ${I} = [e_1, \hdots, e_{k}]$ with $k \le \rk$ be an ordered set of independent elements of $E \setminus f$ such that
\[{(\mathcal{O},e_1,f), (\mathcal{O} / \{e_1\},e_2,f),} \ldots, (\mathcal{O} / \{e_1, \hdots, e_{k-1}\},e_k,f)\] are Euclidean oriented matroid programs.
Let $\mathcal{O}' = \mathcal{O} \cup p$ be the lexicographic extension $\mathcal{O}[e^+_1, \hdots, e^+_k]$.  
Then $(\mathcal{O}',p,f)$ is a Euclidean oriented matroid program.
\end{corollary}

\section{Proof of Theorem \ref{theo:topSweepExists}}\label{sec:five}

\begin{thmn}[\ref{theo:topSweepExists}]
  Let $(\Oscr,g,f)$ be a Euclidean oriented matroid program. Then
  $G_f$ defines a partial order of the cocircuits of $(\Oscr
  \setminus f)$.  Let $f$ be in general position.  Every linear
  extension of that order defines a non-degenerate topological
  sweep and vice versa.
\end{thmn}

\subsection{Sketch of the Proof}

\textcolor{black}{The fact that a topological sweep 
(which is an
 ordering of the cocircuits of $(\Oscr,g,f)$ together with a set of
 parallel extensions going through each cocircuit) 
  always respects
  the partial order defined by $G_f$ follows directly from the
  definition. The other way around is more complicated. We construct a
  sweep by adding parallel extensions going through the cocircuits in
  order of the given linear order.  Using induction, we show that these extensions form a sweep.  We use here the fact that $G_f$ is connected, which we prove in Theorem \ref{theo:graphConnect}.}

\subsection{Topological Sweeping 
(adding results from paper \cite{1})}

We start by considering the paths in the cocircuit graph $G_f$ of an oriented matroid. 
Concerning the connectivity of $G_f$ we cite some definitions and a theorem from \cite{3}.

\begin{definition}[see, \cite{3}, chapter 3]
Let $\mathcal{L}$ be the set of covectors of an oriented matroid $\mathcal{O}$ let $j \in \{+,-\}$ and let $e$ be an element of the groundset of $\mathcal{O}$.
We call $S^j_e = \{X \in \mathcal{L} \mid X_e = j \}$ an {\em open hemisphere} and
$\tilde{S}^j_e = \{X \in \mathcal{L} \mid X_e = j$ or $ 0  \}$ a {\em closed hemisphere} 
of an oriented matroid. We call an intersection of closed hemispheres of an oriented matroid a {\em supercell}.
We call a supercell $\mathcal{S}$ {\em e-positive} if $\mathcal{S} \subseteq \tilde{S}^+_e$ and $\mathcal{S} \cap S^+_e \neq \emptyset $.
A path in the cocircuit graph of an oriented matroid is called {\em e-positive} if all intermediate cocircuits in the path have $e = +$.
\end{definition}

We cite a theorem, see \cite{3}, Lemma 3.3.

\begin{theorem}\label{theo:supercellTheorem}
Let $\mathcal{O}$ be an oriented matroid, and let $e$ be an element of its groundset.
Let $\mathcal{S}$ be an $e$-positive supercell. Then there exists an $e$-positive path between two vertices of the cocircuit graph of $\mathcal{S}$ .
\end{theorem}

\begin{theorem}\label{theo:graphConnect}
Let $(P) = (\mathcal{O},g,f)$ be a feasible and bounded Euclidean oriented matroid program with groundset $E_n \dot{\cup} f \dot{\cup} g$.
Let  $P(E_n)$ be the feasible region of $(P)$, $I$ a covector in $P(E_n)$ and 
let $G_f$ be the graph defined in Definition \ref{def:Graph}. Moreover, assume that $Y_f = +$ for all cocircuits $Y$ of the feasible region. 
\begin{enumerate}
\item The graph $G_f$ is connected.
\item The subgraph $G_{f|I}$ of $G_f$ induced from the cocircuits of $[0,I]$ is connected.
\item The subgraph $G_{f|P(E_n)}$ of $G_f$ induced from the cocircuits of $[0,P(E_n)]$ is connected.
\end{enumerate}
\end{theorem}

\begin{proof}
  ad (i): $G_f$ is the cocircuit graph of an oriented matroid with the
  cocircuits with $g = +$.  It is therefore \textcolor{black}{the
    cocircuit graph of a} supercell. We can apply theorem
  \ref{theo:supercellTheorem} and get a path in $G_f$ between any
  vertices $X,Y$ of $G_f$.  ad (ii): $[0,I]$ is a supercell, see
  \cite{3}, beginning of Chapter 3. If $I_g = 0$, then the subgraph
  $G_{f|I}$ is empty.  If $I_g = +$ it is $g$-positive and two
  vertices in $G_{f|I}$ are connected by a path in $G_{f|I}$ because
  of Theorem \ref{theo:supercellTheorem}.  (iii) follows directly from
  (ii).
\end{proof}

Concerning parallel extensions, we cite as a Lemma a part of the proof
of \cite{2}, Theorem 10.5.5.

\begin{lemma}\label{lem:parallelExtensionEdmondsMandel}
Let $(\mathcal{O},g,f)$ be a Euclidean oriented matroid program let $\mathcal{C}$ the set of its cocircuits and and let $Y^0$ be a cocircuit of \textcolor{black}{$(\mathcal{O},g)$}.
We define a localization $\sigma \colon \mathcal{C} \rightarrow \{+,-,0 \}$ as follows. For any cocircuit $Y \in \mathcal{C}$ holds 
\begin{enumerate} 
\item If $Y_g = +$, then
\[ \sigma(Y) = \begin{cases} + &\text{ if there is a directed path in } G_f \text{ from } Y^0 \text{ to } Y ,   \\
 0 &\text{ if there is an undirected path in } G_f \text{ from } Y^0 \text{ to } Y, \\
 - &\text{ otherwise. }  \end{cases} \]
\item If $Y_g = 0$, then $\sigma(Y) := Y_f$, and
\item If $Y_g = -$, then $\sigma(Y) := -\sigma(-Y)$.
\end{enumerate}
That localization corresponds to an extension of  $(\mathcal{O},g,f)$ with an element $h$ through $Y^0$ parallel to $f$ with respect to $g$. 
\end{lemma}

We add a lemma concerning parallel extensions and the graph $G_f$.

\begin{lemma}[see \cite{1} 10.5.3]\label{lem:parallelExtIdDirectionConform}
Let $(\mathcal{O},g,f)$ be an oriented matroid program, $Y^0 \in (\mathcal{O},g)$ a cocircuit, and $(\mathcal{O}'=\mathcal{O} \cup h,g,f)$ be an extension through $Y^0$ parallel to $f$ wrt. $g$. 
Let $X, Y \in (\mathcal{O},g)$ be two conformal cocircuits such that $X \circ Y$ is an edge. 
Let $X', Y'$ denote the
corresponding cocircuits in $\mathcal{O}'$.
If $\Dir_{g,f}(X,Y) = 0$, we have $X'_h = Y'_h$. If $\Dir_{g,f}(X,Y) = +$ (or $\Dir_{g,f}(X,Y) = -$) we have $X'_h \preceq Y'_h$ (or $Y'_h \preceq X'_h$) and 
 $X'_h =Y'_h = 0$ cannot occur.
\end{lemma}

\begin{proof}
The Lemma follows directly from $\Dir_{g,f}(X,Y) = \Dir_{g,h}(X',Y')$.
\end{proof}

We show a Lemma considering the new cocircuits of the parallel extension.

\begin{lemma}\label{lem:undirectedPathsInParallelExxtensions}
Let $(\mathcal{O},g,f)$ be an oriented matroid program and $(\mathcal{O}',g,f) = (\mathcal{O} \cup h,g,f)$ an extension such that
$h$ is parallel to $f$ wrt. $g$. Let $X,Y$ be two cocircuits in $(\mathcal{O}',g,f)$ with $X_g = Y_g = +$ and $X_h = Y_h = 0$.
Then there is an undirected path in $(\mathcal{O}',g,f)$ from $X$ to $Y$ with only cocircuits with $g = +$.
\end{lemma}

\begin{proof}
Because $X_h = Y_h = 0$, we look at the contraction $(\mathcal{O}' / h,g)$ and see if there is a path. 
$(\mathcal{O}' / h,g)$ is a closed hemisphere with $g$-positive cocircuits, hence a $g$-positive supercell.
We apply Theorem \ref{theo:supercellTheorem} and obtain a $g$-positive path between $X$ and $Y$ in $(\mathcal{O}' / h,g)$
which is then an undirected $g$-positive path in the cocircuit graph of $(\mathcal{O}',g,f)$.
\end{proof}

We add some results from a paper by W. Hochst\"attler \cite{1}. 
We start with a definition.

\begin{definition} \label{def:sweep}[see \cite{1} Definition 8]
Let $(\mathcal{O},g,f)$ be a Euclidean oriented matroid program and let \textcolor{black}{$\mathcal{C}_+$} denote the cocircuits $Y \in (\mathcal{O},g)$ with $Y_g = +$.
A {\em topological sweep of $(\mathcal{O},g,f)$} is a total order $Y^1, \hdots, Y^s$ of \textcolor{black}{$\mathcal{C}_+$}, the {\em sweep order}, together with an extension
$\mathcal{O}' = \mathcal{O} \cup \{h_1, \hdots, h_s\}$ by $s$ elements such that
\begin{enumerate}
\item for all $1 \leq i \leq s$ holds that $h_i$ is parallel to $f$ wrt. $g$
\item for all $1 \leq i \leq s$ there exists $1 \leq i_1 \leq i \leq i_2 $ such that
\[Y'^i_{h_j} = \begin{cases} - & \text{if } j<i_1 \\ 0 & \text{if }i_1 \leq j \leq i_2 \\ + & \text{if } j > i_2 \end{cases}\] 
\end{enumerate}
We call the sweep {\em non-degenerate}, if there is always $i_1 = i_2$.
\end{definition}
\begin{remark}
Recall that in the non-degenerate case all extensions $h_{i_1}, \hdots h_{i_2}$ are parallel to each other.
\end{remark}

We add a definition here.

\begin{definition}
Let $(\mathcal{O},g,f)$ be an oriented matroid program and let $X$ be a cocircuit of $\mathcal{O}$ with $X_g = +$ and $X_f \neq 0$.
If $(\mathcal{O},g,f)$ is Euclidean and we have a topological sweep and a cocircuit $Y^i \in (\Oscr,g)$ with its corresponding extension $h_i$ then we call the oriented matroid program $((\mathcal{O} \cup h_i) \setminus f, g, h_i)$ {\em derived from $(\mathcal{O},g,f)$ by parallel movement of the target function through $Y^i$ respecting the sweep order}.
\end{definition}

The sweep order agrees with the direction of the edges of $G_f$.


\begin{lemma}\label{lem:sweepOrientation}
Let $(\mathcal{O},g,f)$ be a Euclidean oriented matroid program with a topological sweep.
The sweep order is compatible with the orientation of the edges of $G_f$. That means if $X \circ Y$ is an edge in $G_f$ directed from 
$X$ to $Y$, then $X$ precedes $Y$ in the sweep order.
\end{lemma}

\begin{proof}
Let $Y^1, \hdots, Y^r$ be a topological sweep of $(\mathcal{O},g)$. Let $Y^i,Y^j$ be two conformal cocircuits of the sweep with $Y^i \circ Y^j$ being an edge in $G_f$ directed
from $Y^i$ to $Y^j$. In the extension $\mathcal{O} \cup \{h_1, \hdots, h_r\}$ corresponding to the sweep holds $Y^j_{h_j} = 0$.
Strong cocircuit elimination of $g$ holding $h_j$ between $-Y^i$ and $Y^j$ yields a unique cocircuit $Z$ with $Z_g = 0$ and 
$Z_{h_j} = -Y^i_{h_j}$. It holds also $Z_f = +$ because the edge $Y^i \circ Y^j$ is directed from $Y^i$ to $Y^j$.
Because of the parallelity of $h_j$ and $f$ wrt. $g$ we obtain $Z_{h_j} = Z_f = +$ hence $Y^i_{h_j} = - Z_{h_j} = -$.
From Definition  \ref{def:sweep} follows $i < j$ hence $Y^i$ comes in the sweep order before $Y^j$. 
\end{proof}

Because of \cite{1}, Proposition 4, every extension through a cocircuit $X$ parallel to $f$ wrt. $g$ goes through the cocircuits reachable 
from $X$ with an undirected path. 
It follows directly:

\begin{prop}\label{prop:nonDegNoUndirPaths1}
Let $(\mathcal{O},g,f)$ be a Euclidean oriented matroid program. If a non-degenerate topological sweep of $\mathcal{O}$ exists,
the graph $G_f$ contains no undirected edges between cocircuits that have not both $f = 0$.
\end{prop}

But we have to look at what happens with the graph $(G \setminus f)_f$.

\begin{prop}\label{prop:nonDegNoUndirPaths}
Let $(\mathcal{O},g,f)$ be a Euclidean oriented matroid program. If a non-degenerate topological sweep of $\mathcal{O}$ exists,
the graph $(G \setminus f)_f$ contains no undirected edges.
\end{prop}

\begin{proof}
If an edge is directed in $G_f$, it is also directed in $(G \setminus f)_f$ if existing. Edges in $(G \setminus f)_f$ that are not in $G_f$ are also directed, see Definition \ref{def:Graph}.
\end{proof}

\begin{prop}\label{prop:sweepOrientation2}
  Let $(\mathcal{O},g,f)$ be a Euclidean oriented matroid program with
  a topological sweep with $f$ being in general position.  Let $X,Y
  \in (G \setminus f)_f$ cocircuits with $X$ coming before $Y$ in the
  sweep order.  Then there exists no (un-)directed path in $G_f$ from
  $Y$ to $X$.
\end{prop}
\begin{proof}
If there were a path directed from $Y$ to $X$ (or undirected), it would have at least
one edge $(P^i,P^j)$ in $(G \setminus f)_f$ with $i > j$ in sweep order.
But that edge would then have the wrong direction.
\end{proof}

$G_f$ induces a partial order of the cocircuits. Improving the result of Theorem 9 in \cite{1}, we show
that every linear extension of that ordering yields a topological sweep and vice versa.

\begin{thmn}[\ref{theo:topSweepExists}]
  Let $(\Oscr,g,f)$ be a Euclidean oriented matroid program. Then
  $G_f$ defines a partial order of the cocircuits of $(\Oscr
  \setminus f)$.  Let $f$ be in general position.  Every linear
  extension of that order defines a non-degenerate topological
  sweep and vice versa.
\end{thmn}

\begin{proof}
  We have a non-degenerate topological sweep and let $X \neq Y$ be two
  different cocircuits in $(\mathcal{O} \setminus f)$.  We know that
  the sweep respects the order of $G_f$. If $X$ comes before $Y$ in
  the sweep order we set $X \prec Y$ which is an extension of the
  poset.

  On the other hand let $P = P^1, \hdots, P^n$ be a linear extension
  of the partial order induced by $G_f$. We extend
  $(\mathcal{O},g,f)$ in $n$ steps as follows.  Let
  $(\mathcal{O}^{n+1},g,f) = (\mathcal{O},g,f)$ and let $i = n$.
\begin{enumerate}
\item Let $(\mathcal{O}^i,g,f) = (\mathcal{O}^{i+1} \cup h,g,f)$ such that $h$ is an extension 
through $P^i$ parallel to $f$ with respect to $g$ like in Lemma \ref{lem:parallelExtensionEdmondsMandel}.
\item Let $i = i-1$. Goto (i) until $i = 0$.
\end{enumerate}

We show now that $(\mathcal{O}^1,g,f)$ together with $P = P^1, \hdots, P^n$  is the desired sweep.
We use induction to show the following {\em statement A}. 
\[P^j_{h_i} =  \begin{cases} + \text{ if } j > i, \\  - \text{ if } j < i \text{ and } \\ 0 \text{ if } j = i. \end{cases} \]
First let $i = n$. \textcolor{black}{Proposition \ref{prop:sweepOrientation2} yields that} there is no (un)-directed path  in $G_f$ from $P^n$ to a $P^j$ with $j < n$. 
Hence we have $P^j_{h_n} = -$ and $P^n_{h_n} = 0$ in $(\mathcal{O}^n,g,f)$.
Statement A holds for $i=n$.

Let now statement A hold for all $i+1, \hdots, n$. Let $P^j$ be a cocircuit with $j > i$.
Then because $G_f$ is connected there is a path from $P^i$ to $P^j$ in $(\mathcal{O},g,f)$. Let $P^i \circ Y$ be the first edge of that path.
We have $Y_g = +$. The edge is directed because $f$ is in general position. 
Cocircuit elimination of $g$ between $-P^i$ and $Y$ yields a cocircuit $Z'$ with $Z'_f \neq 0$.
If $Z'_f = -$ let $Z = -Z'$ otherwise let $Z = Z'$. 
In $(\mathcal{O}^{i+1},g,f)$ we have $Z_{h_j} = Z_f = +$ because of the parallelity of $f$ and $h_j$.
Also we have $P^i_{h_j} = -$ because statement A holds for $(\mathcal{O}^{i+1},g,f)$.  
$P^i \circ Z$ is still an edge in $(\mathcal{O}^{i+1},g,f)$.
Cocircuit elimination of $h_j$ between $P^i$ and $Z$ yields a cocircuit $X$ and a directed path 
in the cocircuit graph $G^{i+1}_f$ of $(\mathcal{O}^{i+1},g,f)$ from $P^i$ to $X$.
Because $X_{h_j} = P^j_{h_j} = 0$ Lemma \ref{lem:undirectedPathsInParallelExxtensions}
yields an undirected path in $G^{i+1}_f$ from $X$ to $P^j$ hence we obtain a directed path
from $P^i$ to $P^j$ in $G^{i+1}_f$ and we have $P^j_{h_i} = +$ in $(\mathcal{O}^i,g,f)$.

Now let $P^j$ be a cocircuit with $j < i$ in
$(\mathcal{O}^{i+1},g,f)$.  Because statement A holds in
$(\mathcal{O}^{i+1},g,f)$ we have $P^j_{h_{i+1}} = P^i_{h_{i+1}} = -$.
A directed path from $P^i$ to $P^j$
in $G^{i+1}_f$ can not have cocircuits with $h_{i+1} \in \{0,+\}$
because otherwise we would obtain two edges in the path of opposite
direction \textcolor{black}{because for an edge $X \circ Y$ with
  $X_{h_{i+1}} = -$ and $Y_{h_{i+1}} = 0$ holds $\Dir_{g,f}(X,Y) = +$}
.  Hence all cocircuits in the path must have $h_j = -$ for all $j >
i$.  But that means the cocircuits are old, exist already in
$(\mathcal{O},g,f)$ and form a directed path in $G_f$ (\see
Proposition \ref{prop:EdgeStaysEdge} and Lemma
\ref{lem:DirPresOfOutsideEdges}) But that path would
  at least contain one edge $(P^k, P^l)$ in $G_f$ with $k>l$ directed
  from $P^k$ to $P^l$ violating the fact that $P$ is a linear
  extension of the partial order of $G_f$.  We obtain $P^j_{h_i} = -$
for all $j < i$ in $(\mathcal{O}^i,g,f)$.

Finally $P^i_{h_i} = 0$ in$(\mathcal{O}^i,g,f)$ follows directly from the construction.
Statement A holds for  $(\mathcal{O}^i,g,f)$ and hence for $(\mathcal{O}^1,g,f)$ which proves everything.
\end{proof}

\section{Proof of Theorem \ref{theo:FourthMainTheoem}}\label{sec:four}

\begin{thmn}[\ref{theo:FourthMainTheoem}]
 Let $\mathcal{O}$ be a Euclidean oriented matroid with $f \in \mathcal{O}$ lying in general position.
 Then each tope cell of $\mathcal{O} \setminus f$ is vertex-shellable.
\end{thmn}

\subsection{Sketch of Proof}

It is enough to show the theorem for the feasible region $P$ of a
feasible and bounded oriented matroid program $(\mathcal{O}, g, f)$
because each tope cell of $\mathcal{O} \setminus f$ can (maybe after
reorienting some elements) be seen like that.  We take a sweep $Y^1,
\hdots, Y^n$ of the cocircuits of $P$ respecting the order of the
cocircuit graph of $(\mathcal{O}, g, f)$ and prove that it is a
vertex-shelling.  We show in Subsection
\ref{Subsection:SweepOrderAndGraph} that for each cocircuit $Y^j$ the
intersection of the interval $[Y^j,P]$ with the union of the intervals
$[Y^i,P]$ for all $i < j$ equals a union of some intervals $[Y^k \circ
Y^j,P]$ with $k \in J_j$ and $J_j \subseteq \{1, \hdots, j-1 \}$. That
is the non-recursive part of the definition of a vertex-shelling.  In
Subsection \ref{Subsection:theRecursiveAtomOrdering} we handle the
recursive part of the definition. We show that there is again a
vertex-shelling of the interval $[Y^j, P]$ such that in the ordering
of its atoms (which are edges in the original oriented matroid), the
edges $Y^k \circ Y^j$ with $k \in J_j$ come first.  To achieve that
recall that in the sweep order of a Euclidean oriented matroid program
$(\mathcal{O}, g, f)$ the cocircuits with $f = -$ come first.  Moving
the target function $f$ over the cocircuit $Y_j$ yields the oriented
matroid program $(\mathcal{O}' = \mathcal{O} \cup f' \setminus f, g,
f')$.  Applying induction on the rank we may assume that in that
program the interval $[Y^j, P]$ has a vertex-shelling where the
cocircuits with $f' = -$ come first.  But these are exactly the edges
$Y^j \circ Y^k$ with $k \in J_j$ in the original ordering.  To get the
inductive assumption fulfilled we have to solve the problem that in
the interval $[Y^j, P]$ of $\mathcal{O}'$ the element $g$ does not
exist anymore.  We have to substitute  $g$ with another element
$e$ that exists in the minor. Moving the target function $f$ in
$(\mathcal{O}, g, f)$ preserves Euclideaness in $(\mathcal{O}', g,
f')$, but it does not necessarily preserve the Euclideaness of another
program $(\mathcal{O}', e, f')$.  We handle that problem in Subsection
\ref{Subsection:preservingEuclideanessInMinors}.  But 
substituting the element $g$ with another element still is not quite enough
to fulfill the inductive assumption because we could maybe lose the
boundedness of the oriented matroid program in the minor.  We solve
that problem in the proof of the main theorem in Subsection
\ref{Subsection:theRecursiveAtomOrdering} using a lexicographic
extension preserving Euclideaness.

\subsection{Preservation of Euclideaness in Minors from Parallel Movements of the Target Function}\label{Subsection:preservingEuclideanessInMinors}

For the definition of Euclideaness of oriented matroids and parallel single-element extensions, see \cite{2}, 10.5.1, 10.5.2.
Concerning the Dir-function, we get immediately from the definition of parallel extensions.

\begin{prop}\label{prop:parallelExtensionsDoNotChangeDirection}
Let $(\mathcal{O},g,f)$ be a Euclidean oriented matroid program.
Let $\mathcal{O} \cup f'$ be a single-element extension parallel to $f$ wrt. $g$. Then it holds
\[ \Dir_{g,f}(X,Y) = \Dir_{g,f'}(X,Y) \text{ for all cocircuits } X,Y \in (\mathcal{O},g). \]
\end{prop}

First, we add an important definition.

\begin{definition}
Let $(\mathcal{O},g,f)$ be an oriented matroid program and let $X$ be a cocircuit of $\mathcal{O}$ with $X_g = +$ and $X_f \neq 0$. Let $\mathcal{O} \cup f'$ be an extension of $\mathcal{O}$ parallel to $f$ wrt. $g$ satisfying $\sigma(X)=0$.
Then we call the oriented matroid $((\mathcal{O} \cup f') \setminus f)$ {\em derived from $(\mathcal{O},g,f)$ by parallel movement of the target function (through $X$)}.
\end{definition}

The next lemma follows directly.

\begin{lemma}\label{lem:parVersch}
Let $(\mathcal{O},g,f)$ be an oriented matroid program with graph $G_f$. Let $(\mathcal{O}' = \mathcal{O} \cup h)$ be derived from $(\mathcal{O},g,f)$ by parallel movement of the target function. Then it holds for the graph $G'_h$ of $(\mathcal{O}',g,h)$ that $(G \setminus f)_f = (G' \setminus h)_h$.
If $(\mathcal{O},g,f)$ is Euclidean, feasible and bounded, so is $(\mathcal{O},g,h)$.
\end{lemma}

\begin{proof}
The identity of the graphs follows from the identity of the graphs $G''_f = G''_h$ of the oriented matroid programs $(\mathcal{O} \cup \{f,h\},g,f)$ and $(\mathcal{O} \cup \{f,h\},g,h)$. 
The feasible region of $(\mathcal{O},g,f)$ and $(\mathcal{O},g,h)$ stays the same, Euclideaness and boundedness are also clear.
\end{proof}

We add a small remark.
\begin{remark}\label{prop:gZeroCocircuitsIdenticalMovingTargetFunction}
If $(\mathcal{O},g,f)$ is an oriented matroid program and $(\mathcal{O}' = \mathcal{O} \cup f' \setminus f)$ is derived from $(\mathcal{O},g,f)$ by parallel movement of the target function, then the cocircuits with $g = 0$ are identical in $\mathcal{O}$ and in $\mathcal{O}'$ because it holds $\mathcal{O} / g = \mathcal{O}'/ g$.
\end{remark}

%

Now we show that in a special situation Euclideaness is preserved in minors of Euclidean oriented matroid programs derived by parallel movement of the target function.
We define the following objects.

\begin{definition}
 Let $\mathcal{O}^1$ be a Euclidean oriented matroid of rank $r$. 
  Let $f \neq g \in \mathcal{O}^1$ be  elements in general position.
 Let $X$ be a cocircuit of $\mathcal{O}^1$ with $X_g = +$ and $X_f \neq 0$.
 \begin{enumerate}
 \item Let $\mathcal{O}^2 = (\mathcal{O}^1 \cup f') \setminus f$ derived from $\mathcal{O}^1$ by parallel movement of the target function 
 through $X$ (hence $X_{f'} = 0$ in $\mathcal{O}^2$). 
  \item Let $\mathcal{O}^3 = (\mathcal{O}^2 \setminus (\supp(X)\setminus \{ g\})$ hence $f' \in \mathcal{O}^3$ .  
  \item Let $ \mathcal{O}^4 = (\mathcal{O}^1 \setminus (\supp(X) \setminus \{f,g\})$.
  \item Let $\mathcal{O}^5 = (\mathcal{O}^4 \cup f^{''}) \setminus f$ be derived from $(\mathcal{O}^4,g,f)$ by parallel movement of the target function through $X$.
\end{enumerate}

\end{definition}

We consider the following diagram:
\begin{center}
\includegraphics[scale=0.7]{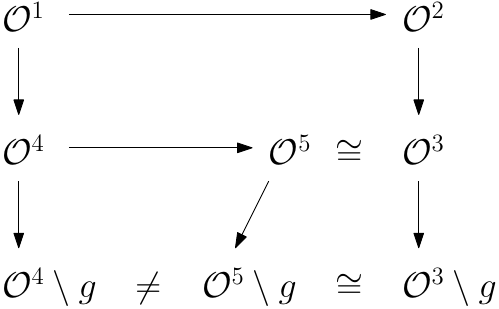} 
\end{center}

Considering the ranks of these objects, we get:

\begin{prop}
All objects in the diagram have rank $r$ except $\mathcal{O}^5 \setminus g$ and  $\mathcal{O}^3 \setminus g$
which have rank $r-1$.
\end{prop}
 
\begin{proof}
  All oriented matroids of the diagram contain $z(X)$ hence they have
  at least rank $r-1$ but all objects except $\mathcal{O}^3 \setminus g$ contain elements $\notin z(X)$
  hence have rank $\geq r$ but because $\mathcal{O}$ has rank $r$ and
  a parallel extension does not increase the rank they have all rank equal to $r$.
\end{proof} 
 
Because $f$ is in general position in $\mathcal{O}^1$ we obtain:

\begin{prop}
\
\begin{enumerate}
\item $f$ is in general position in $\mathcal{O}^4$. 
\item  $f'$ is in $X$ in general position in $\mathcal{O}^2$ hence also in $\mathcal{O}^3$. 
\end{enumerate}
\end{prop}

\begin{proof}
Only (ii) needs a proof. If $f'$ would not lie in general position in $X$ in $\mathcal{O}^2$, then it would lie in an edge $F$ with $f' \in \cl(z(F) \setminus f')$. 
We do not have $g \in z(F)$ because otherwise we would have $X_g = 0$ which is impossible.
But then there is a cocircuit $Y$ with $z(Y) \supseteq z(F) \cup g$. We have $Y_{f'} = 0$ hence $Y_f = 0$ hence $f$ is not in general position.
\end{proof}

We start by analyzing the smaller objects:

\begin{prop} \label{prop:12}
\
\begin{enumerate} 
\item In
$\mathcal{O}^5$ resp. 
$\mathcal{O}^3$ we have $X_{f''} = X_{f'} = 0$ and $X_g \neq 0$. All the other cocircuits $Y$ have  $Y_g = 0$ and $Y_{f''} = Y_{f'}$.
\item For a cocircuit $Y$ in $\mathcal{O}^4$ with $Y_f \neq 0$ holds $Y = \pm X$ or $Y_g = 0$.
\end{enumerate}
\end{prop}

\begin{proof}
  If $Y \neq \pm X$ in $\mathcal{O}^3$ (or $\mathcal{O}^5$) we have $z(Y) \neq z(X)$ hence there is a $h \in z(Y)$ with $h
  \notin z(X)$. Thus, $h = g$ and we obtain $Y_g = 0$.
  This and Remark \ref{prop:gZeroCocircuitsIdenticalMovingTargetFunction} prove (i). Statement (ii) works analogously using the fact
  that $X_f \neq 0$ in $\mathcal{O}^4$.
\end{proof}

We obtain directly.
\begin{corollary}\label{cor:55} It holds $\mathcal{O}^3 \cong \mathcal{O}^5$ hence also $\mathcal{O}^3 \setminus g \cong \mathcal{O}^5 \setminus g$.
The isomorphism maps $f'$ to $f^{''}$ and is the identity for the rest.
\end{corollary}
\begin{proof}
It holds $\mathcal{O}^3 \setminus f' = \mathcal{O}^5 \setminus f''$ and the first statement follows from Proposition \ref{prop:12} (i).
The second statement is immediate.
\end{proof}
Because of that Corollary we will identify $\mathcal{O}^5$ and $\mathcal{O}^3$ in the following.

\begin{prop} \label{prop:die13er}
\[ \mathcal{O}^4 \setminus f  = \mathcal{O}^3 \setminus f' \]
\end{prop}

\begin{proof} This follows from
\[ \mathcal{O}^4 \setminus f = (\mathcal{O}^1 \setminus \supp(X)) \cup g = \mathcal{O}^2 \setminus \{supp(X) \cup f'\} \cup g = \mathcal{O}^3 \setminus f'. \qedhere \]
\end{proof}

\begin{lemma}\label{lem:66}
  The cocircuits of $\mathcal{O}^4$ and $\mathcal{O}^3$ with $g = 0$
  are identical.  The cocircuits with $f,f' \neq 0$ are also identical
  except $X$ with $X_{f'} = 0$ in $\mathcal{O}^3$ and $X_f \neq 0$ in
  $\mathcal{O}^4$.
\end{lemma}

\begin{proof}
The first statement is Remark \ref{prop:gZeroCocircuitsIdenticalMovingTargetFunction}.
The second statement follows from Proposition \ref{prop:12}.
\end{proof}

We start to compare the Euclideaness of the objects.

\begin{lemma}\label{lem:o4euclo5eucl}
$(\mathcal{O}^3,g,e)$ is Euclidean for all $e \in \mathcal{O}^3 \setminus g$.
\end{lemma}

\begin{proof}
Let $Y$ be a cocircuit in $\mathcal{O}^3$ with $Y_g = +$. Then we have $z(Y) \subseteq z(X)$ hence $Y = X$,
we cannot have directed cycles in  $(\mathcal{O}^3,g,e)$. 
\end{proof}

\begin{lemma}\label{lem:o4euclo5eucl1}
If $(\mathcal{O}^4,e,f)$ is Euclidean for an $e \in \mathcal{O}^4 \setminus \{f,g\}$ then is also $(\mathcal{O}^5,e,f'')$.
\end{lemma}

\begin{proof}
We assume there is a directed cycle in the graph $G_{f''}$ of $(\mathcal{O}^5,e,f'')$.
Then we know from Proposition \ref{prop:PathGraphWellKnownFacts} (iv) that a directed cycle can only have cocircuits with $f'' \neq 0$ hence
we know from Proposition  \ref{prop:12} that it has to hold $g = 0$ for such cocircuits.
But from Lemma \ref{lem:66} we know that these cocircuits $Y$ correspond to cocircuits $Y'$ in $\mathcal{O}^4$
with $Y_{f''} = Y'_f$. Edges $Y^1 \circ Y^2$ with $\sep(Y^1,Y^2) = \emptyset$ and $Y^1_g = Y^2_g = 0$ in $\mathcal{O}^5$
correspond to edges $Y'^1 \circ Y'^2$ with $\sep(Y'^1,Y'^2) = \emptyset$ and $Y'^1_g = Y'^2_g = 0$ in $\mathcal{O}^4$ and because cocircuit elimination of $e$ 
between $-Y^1$ and $Y^2$ yields a cocircuit $Z$ with $Z_g = 0$ in $\mathcal{O}^5$ cocircuit elimination of $e$
between $Y'^1$ and $Y'^2$ yields the corresponding cocircuit $Z'$ with $Z'_g = 0$ in $\mathcal{O}^4$ and $Z'_f = Z_{f''}$. The edge $Y'^1 \circ Y'^2$ has the same 
direction like the edge $Y^1 \circ Y^2$
hence we obtain a directed cycle in $\mathcal{O}^4$ which proves everything.
\end{proof}

\begin{lemma}\label{lem:o4euclo5eucl21}
If $(\mathcal{O}^4,e_1,e_2)$ is Euclidean for $e_1 \neq e_2 \in \mathcal{O}^4 \setminus \{f,g\}$ then is also $(\mathcal{O}^5,e_1,e_2)$.
\end{lemma}

\begin{proof}
We assume we have a directed cycle in $(\mathcal{O}^5,e_1,e_2)$.
$X$ is not involved in the cycle since we have $e_1 \notin \supp(X)$ hence  $X_{e_1} = 0$. 
Because of Proposition \ref{prop:12} (i) we have only cocircuits in the cycle with $g=0$ and have a directed cycle in $(\mathcal{O}^4,e_1,e_2)$ like in the lemma before.
\end{proof}

Now we put everything together, we repeat again the assumptions.

\begin{theorem}\label{:theo:PreservationOfEuclideaness}
\begin{enumerate}
 \item Let $\mathcal{O}^1$ be an oriented matroid of rank $r$. 
 \item Let $f \neq g \in \mathcal{O}^1$ be elements in general position.
  \item Let $X$ be a cocircuit of $\mathcal{O}^1$ with $X_g = +$ and $X_f \neq 0$. 
 \item Let $\mathcal{O}^2 = (\mathcal{O} \cup f') \setminus f$ derived from $(\mathcal{O}^1,g,f)$ by parallel movement of the target function 
 through $X$ (hence $X_{f'} = 0$ in $\mathcal{O}^2$).
  \item Let $\mathcal{O}^3 = (\mathcal{O}^2 \setminus \supp(X)) \cup \{g\}$.  
\end{enumerate}
If $\mathcal{O}^1$ is Euclidean then $\mathcal{O}^3$ is a Euclidean oriented matroid, too. If  $(\mathcal{O}^1,e,f)$ is a Euclidean oriented matroid program
for an $e \in $ \textcolor{black}{$z(X)$} then $(\mathcal{O}^3 \setminus g,e,f')$ is a Euclidean oriented matroid program, too.
\end{theorem}

\begin{proof}
$(\mathcal{O}^1,e,f)$ resp. $(\mathcal{O}^1,e,e_2)$ is a Euclidean oriented matroid program 
for all $e \neq e_2 \in \mathcal{O}^1 \setminus f$. 
Let $e \neq e_2  \in \mathcal{O}^3 \setminus \{g,f'\}$.
We know that $(\mathcal{O}^4,e,f)$ resp.$(\mathcal{O}^4,e,e_2)$ is Euclidean as a minor of $(\mathcal{O}^1,e,f)$ resp. $(\mathcal{O}^1,e,e_2)$.
From Corollary \ref{cor:55}, Lemmas \ref{lem:o4euclo5eucl}, \ref{lem:o4euclo5eucl1} resp.\ \ref{lem:o4euclo5eucl21} 
follows the Euclideaness of $(\mathcal{O}^5,g,e)$, $(\mathcal{O}^5,e,f'')$ resp. $(\mathcal{O}^5,e,e_2)$.  
The Euclideaness of $(\mathcal{O}^5,g,f'')$ follows from the Euclideaness of $(\mathcal{O}^4,g,f)$ because of Lemma \ref{lem:parVersch}.
Hence $\mathcal{O}^5$ is a Euclidean oriented matroid and because of the isomorphism of Corollary \ref{cor:55} that holds also for  $\mathcal{O}^3$.
The last statement follows analogously using Lemma \ref{lem:o4euclo5eucl1}, Corollary \ref{cor:55} and the fact that minors of Euclidean matroid programs stay Euclidean.
\end{proof}

\subsection{The Sweep Orders and the Graph $G_f$}\label{Subsection:SweepOrderAndGraph}

\begin{lemma}\label{lem:optimalSolution}
Let $(\mathcal{O},g,f)$ a feasible, bounded Euclidean oriented matroid program with a topological sweep. 
Let $Y^0, \hdots, Y^s$ be the cocircuits of the feasible region $P(E_n)$ of $\mathcal{O} \setminus f$  in sweep order.
Moreover, we assume $Y^0_f = -$ (resp. $Y^s_f = +$). 
\begin{enumerate}
\item Then $Y^0$ is a minimal resp. $Y^s$ is an maximal solution.
\item Let $X,Y$ be two cocircuits of $P(E_n)$, which both are minimal (resp.\ maximal) solutions.
Then there is an undirected path in $G_{f|P(E_n)}$ from $X$ to $Y$.
\item If the sweep is non-degenerate $Y^0$ (resp. $Y^s$) is the unique minimal (maximal) solution of the program. 
\end{enumerate}
\end{lemma}

\begin{proof}
  ad (i): There exist no $Y^i \neq Y^s$ with a conformal edge
  directed from $Y^s$ to $Y^i$ otherwise Lemma
  \ref{lem:sweepOrientation} would imply  $i > s$.  Also there does not exist a cocircuit $Y$ with
  $f=0$ with a conformal edge directed from $Y^s$ to $Y$.  From
  \cite{2}, Lemma 10.1.15 and the boundedness of $(\mathcal{O},g,f)$
  follows that $Y^s$ is a maximal solution and analogously that $Y^0$
  is a minimal solution.

ad (ii): If $X$ and $Y$ are both maximal (minimal) solutions, then $X \circ Y$ and all cocircuits contained in the interval $I = [0,X \circ Y]$ are maximal (minimal) solutions,
see \cite{6}, (2:26) for a proof of that.
The graph $G_{f|I}$ is also connected, which means that there is a path from $X$ to $Y$ with only edges and cocircuits of $[0,I]$.
Since all cocircuits of $G_{f|I}$ are maximal (minimal) solutions, the edges between them must be undirected, that means there exists an undirected path from $X$ to $Y$ in $G_{f|I}$ hence in $G _f$.
(iii) follows directly from (ii). 
\end{proof}

We need here a small lemma for the preservation of Euclideaness in contractions.

\begin{lemma}\label{lem:ContractionStaysBounded}
Let  $P = (\mathcal{O},g,f)$ be a feasible and bounded Euclidean oriented matroid program.
 Let $S \subseteq E \setminus \{f,g\}$.
Then  $P' = (\mathcal{O} / S,g,f)$ is a bounded Euclidean oriented matroid program.
\end{lemma}

\begin{proof}
Minors of Euclidean matroids are Euclidean, see \cite{2}, Corollary 10.5.6.
$(P')$ is also bounded because if there would exist a cocircuit $Y'$ in the feasible region of $(P')$ with $Y'_g = 0$, 
there would exist a corresponding cocircuit $Y$ in $(P)$ with $Y_e = Y'_e$ for all $e \in E \setminus S$ and $Y_e = 0$ for $e \in S$. Then $Y$ would lie in $P(E_n)$ with $Y_g = 0$ contradicting the boundedness of $(\mathcal{O},g,f)$. 
\end{proof}

Now we obtain our first main result concerning the recursive atom-ordering.

\begin{theorem}\label{theo:firstMainResult}
  Let \textcolor{black}{$P = $} $(\mathcal{O},g,f)$ a feasible,
  bounded Euclidean oriented matroid program with a non-degenerate
  topological sweep.  Let $Y^0, \hdots, Y^s$ be the the cocircuits of
  the feasible region $P(E_n)$ of $\mathcal{O} \setminus f$ in sweep
  order.
\begin{enumerate}
\item For each $j > 0$ there exists a directed path in the feasible
  region of $\mathcal{O} \setminus f$ from $Y^0$ to $Y^j$.
\item For all $j > 0$ there exists $i < j$ and a directed edge 
from $Y^i$ to $Y^j$. 
\item For all $j,k$ with $0 \leq k < j \leq s$ there exists a directed edge 
from a $Y^i$ to $Y^j$ with $i < j$ such that $Y^i \circ Y^j \preceq Y^k \circ Y^j$. 
\item For all $0 \leq j \leq s$ let 
\[ Q_j = \left \{ Y^i \circ Y^j \text{ for all } i < j \text{ where } Y^i \circ Y^j \text{ is an edge with } \sep(Y^i,Y^j) = \emptyset \text{ in }  \mathcal{O} \setminus f  \right\} \]
Then it holds (also if we substitute $1$ by $P(E_n)$)
\[ [Y^j,1] \cap \left(\bigcup_{i<j}[Y^i,1]\right) = \bigcup_{Y \in Q_j}[Y,1]\] 
\item If we have $Y^j_f = 0$ and $Y^i_f = -$ for all $i < j$ then it holds  
\[ Q_j = \{ Z \text{ such that } Z_f = - \text{ and }Z\text{ is an atom in }[Y^j,P(E_n)]\}. \]
\end{enumerate}
\end{theorem}

\begin{proof}
We will assume throughout the proof that $Y^j$ is a cocircuit with $0<j$. 
Because parallel movement of the target function does neither change the cocircuits $Y^0, \hdots, Y^s$ nor their sweep order 
nor the feasibility or boundedness of $(\mathcal{O},g,f)$,
we may assume that $Y^j_f = 0$, $Y^i_f = -$ for $i < j$ and $Y^i_f = +$ for $i > j$.  
$Y^0$ is the unique minimal solution of the program $(\mathcal{O},g,f)$. Note, that we assume that $Y^j$ is a degenerate vertex in $\mathcal{O}$ and ``still is'' a vertex in $\mathcal{O}\setminus f$.

ad (i): Let $\mathcal{D}$ be the set of cocircuits of $P(E_n)$ from
which there exist a directed or undirected path to $Y^j$ (note, that
$Y^j \in \mathcal{D}$). Let $Y^i \in \mathcal{D}$ be the element which
is minimal in the sweep order.  By Lemma
\ref{lem:sweepOrientation} there are no cocircuits $X$ in the feasible
region of $(\mathcal{O},g,f)$ with $X \circ Y^i$ being an edge in
$G_f$ directed from $X$ to $Y^i$. By \cite{2}, Lemma 10.1.15,
and  the boundedness of $(\mathcal{O}',g,f')$ it follows that $Y^i$
is a minimal solution in $(\mathcal{O},g,f)$ hence $Y^i = Y^0$.  There
is a directed path in $(\mathcal{O},g,f)$ from $Y^0$ to $Y^j$. Because
all cocircuits in the path except $Y^j$ have $f = -$, they ``are'' also
cocircuits in $\mathcal{O} \setminus f$. Hence we have  a
directed path already there.

ad (ii): 
For each $Y^j$ there exist a directed path from $Y^0$ to $Y^j$ in $\mathcal{O} \setminus f$. The last edge in this path is the desired edge (it is directed because we assume the sweep to be non-degenerate). 

ad (iii): See also \cite{2}, Lemma 4.1.8 and the remarks after the lemma. We use contraction and consider the program 
\[(P') = (\mathcal{O} / z(Y^k \circ Y^j),g,f).\]
Note that $f,g \in \supp(Y^k \circ Y^j)$.
Its face lattice is isomorphic to the intervall $[0,Y^k \circ Y^j]$. 

$(P')$ is an oriented matroid program which is feasible, because
$z(Y^j)$ and $z(Y^k)$ are supersets of $z(Y^k \circ Y^j)$, hence $Y^j$
and $Y^k$ are also cocircuits of $(P')$ lying in the feasible region
of $(P')$.  It is also bounded and Euclidean because of Lemma
\ref{lem:ContractionStaysBounded} hence there is a topological sweep
of the cocircuits of $(P')$.  Because the edges between the cocircuits
in the graph $G_{f|P'}$ of $(P')$ are oriented like in the graph $G_f$
in $(P')$ and because of Theorem \ref{theo:topSweepExists} we can choose
the sweep order of the cocircuits of $(P')$ as a restriction of the
sweep order of $(P)$.  Hence (ii) holds for $(P')$ and there is an
edge from an $Y^i$ to $Y^j$ in $(P')$ with $i < j$. But this is also
an edge in $(P)$ lying in $[0,Y^k \circ Y^j]$ which proves everything.

ad (iv): We always have
\[ [Y^j,1] \cap (\bigcup_{i<j}[Y^i,1]) = \bigcup_{i<j}[Y^i \circ Y^j,1] \text{ for all } 0 \leq j \leq r.\] 
If $Y^k \circ Y^j$ is not an edge for $k<j$ then because of (iii) we have another $i<j$ such that $Y^i\circ Y^j$ is an edge and $Y^i,Y^j$ conformal (hence $Y^i \circ Y^j \in Q_j$) and 
$[Y^k \circ Y^j,1] \subset [Y^i \circ Y^j, 1]$. The equation from (iv) follows. The substitution of $1$ by $P(E_n)$ works analogously.

ad (v): Let $Y^j_f = 0$ and $Y^i_f = -$ for all $i < j$.  Let $Z \in
Q_j$. Then since $Z = Y^i \circ Y^j$ is an edge for an $i<j$ it is an
atom in $[Y^j,P(E_n)]$ and
$Z_f = -$.  On the other hand, if $Z$ is an atom of $[Y^j,P(E_n)]$ and $Z_f
= -$, then $Z = Y^i \circ Y^j$ for a $i \neq j$ with $Y^i \circ Y^j$
being an edge and $Y^i, Y^j$ being conformal cocircuits.  Clearly, we
must have $i<j$.
\end{proof}

\subsection{The Recursive Atom-ordering}\label{Subsection:theRecursiveAtomOrdering}

So far we have proven 'first level' in the shellability of the feasible region of the
program in  Theorem \ref{theo:firstMainResult}. For shellability, we have to show that we can proceed recursively this way.

We start with an important Lemma.
  
\begin{lemma}\label{lem:EdgesAreCocircuits}
Let $(\mathcal{O}^1,g,f)$ be a feasible, bounded Euclidean oriented matroid program of rank $r > 2$.
Let $X$ be a cocircuit of the feasible region $P(E_n)$ of $\mathcal{O}^1$ 
with $X_f = 0$ and $X \setminus f$ being a cocircuit in $\mathcal{O}^1 \setminus f$ and let $\mathcal{O}^3 = \mathcal{O}^1 \setminus \supp(X)$.
 Let 
\[f \colon \mathcal{O}^3 \rightarrow [X,1] : Y \mapsto Y' \text{ such that } Y'_e = \begin{cases} Y_e \text{ if } e \in  z(X) \\ X_e \text{ if } e \in \supp(X) \end{cases}\]
and let $f^{-1}(Y') = Y' \setminus \supp(X)$. Then $f$ is an isomorphism from $\mathcal{O}^3 \setminus f$ to $[X,1]$ in $\mathcal{O}^1 \setminus f$ mapping
 the feasible region $P'$ of $\mathcal{O}^3 \setminus f$  to $[X,P(E_n)]$ in $\mathcal{O}^1 \setminus f$.
\end{lemma}  

\begin{proof}
See also \cite{2}, Proposition 4.1.9.
It is clear that $f^{-1}(Y') \in \mathcal{O}^3$ if $Y' \in  [X,1]$.
Let now $Y \in \mathcal{O}^3$. Then we have $Y = Z \setminus \supp(X)$ for a covector $Z \in \mathcal{O}^1$.
Then we have $X \circ Z = f(Y) \in [X,1]$. Hence $f$ is a bijection of covectors.
Also it holds $Y^1 \prec Y^2$ iff $f(Y^1) \prec f(Y^2)$ hence $f$ is an isomorphism.
If now $Y$ is in the feasible region of $\mathcal{O}^3$ then (because $X$ is in the feasible region of $\mathcal{O}^1$)
also $f(Y)$ lies  in the feasible region of $\mathcal{O}^1$ and vice versa. 
\end{proof}

We recall the definition of a recursive atom-ordering.
 
 \begin{definition}[see \cite{2}, Definition 4.7.17]\label{def:recAtomOrdering}
 Let $P$ be a graded poset. A linear ordering $x_1,x_2, \hdots, x_s$ of its atoms is a {\em recursive atom-ordering} if either $length(P) \leq 2$,
 or if $length(P) > 2$ and for all $1 \leq j \leq s$ there is a distinguished subset $Q_j$ of the atoms of $[x_j,1]$ such that
 \begin{enumerate}
 \item $[x_j,1] \cap (\bigcup_{i<j}[x_i,1]) = \bigcup_{y \in Q_j}[y,1]$, and
 \item $[x_j, 1]$ has a recursive atom-ordering in which the elements of $Q_j$ come first.
 \end{enumerate}
 \end{definition}

Now we are prepared for our main theorem:

 \begin{theorem}\label{theo:VertexShellingExist}
\begin{enumerate}
 \item Let $\mathcal{O}^1$ be a Euclidean oriented matroid of rank $\rk$. 
   \item Let $f \neq g \in \mathcal{O}^1$ be elements in general position and not coloops.
 \item Let $(\mathcal{O}^1,g,f)$ be a feasible, bounded Euclidean oriented matroid program.
\item  $(\mathcal{O}^1,g,f)$ has a non-degenerate topological sweep because of Theorem \ref{theo:topSweepExists}.
 \item Let $Y^0,\hdots, Y^s$ be the cocircuits of the feasible region $P(E_n)$ in $\mathcal{O}^1 \setminus f$
 in sweep order.
\end{enumerate}
Then the sweep-order $Y^0,\hdots, Y^s$  induces a recursive atom-ordering of the face lattice $\mathcal{F}$ of $P(E_n)$ in $\mathcal{O}^1 \setminus f$. 
\end{theorem}
  
\begin{proof}
We prove the theorem by using induction over the rank $\rk(\mathcal{F})$.
If $\rk(\mathcal{F}) \leq 2$ there is nothing to prove. Now let the claim be true for $\rk(\mathcal{F}) \leq \rk-1$ and we assume $\rk(\mathcal{F}) = \rk$. 

We consider a cocircuit $Y^j \in P(E_n)$ with $j \leq s$ and let
 $(\mathcal{O}^2,g,f')$ be the parallel movement of the target function through $Y^j$ respecting the sweep order. 
 $(\mathcal{O}^2,g,f')$ is also a Euclidean oriented matroid program and $\mathcal{O}^2 \setminus f'$ has the same feasible region. 
It holds $Y^i_{f'} = -$ for $i < j$ and $Y^i_{f'} = +$ for $i > j$ in  $(\mathcal{O}^2,g,f')$.
We know from Theorem \ref{theo:firstMainResult} that the $Y^0, \hdots, Y^s$ and the set 
\[ Q_j = \{ Y^i \circ Y^j \text{ for all } i < j \text{ where }  Y^i \circ Y^j \text{ is an edge with } \sep(Y^i \circ Y^j) = \emptyset \text{ in } \mathcal{O} \setminus f   \} \]
 fullfill the part (i) of Definition \ref{def:recAtomOrdering}.
Theorem \ref{theo:firstMainResult} yields also
\begin{equation}\label{equ:atomsQj}
Q_j = \{Z \text{ is an atom in }[Y^j, P(E_n)] \text{ such that } Z_{f'} = -\}. 
\end{equation}

Now let  $\mathcal{O}^3 = (\mathcal{O}^2 \setminus \supp(Y_j))$. We have $\mathcal{O}^3 \setminus f' = $ \textcolor{black}{$ (\mathcal{O}^1  \setminus \supp(Y_j)) \setminus f$} and $\rk(\mathcal{O}^3) = \rk-1$.
 Let $S' = [e_1, \hdots, e_{r-1}]$ be an ordered set of independent elements of $\mathcal{O}^3 \setminus f'$
and let $(\mathcal{O}^3 \cup p)$ be the lexicographic extension $\mathcal{O}^3[e^+_1, \hdots, e^+_{r-1}]$.
\begin{enumerate}
\item $\mathcal{O}^3$ is Euclidean  because of Theorem \ref{:theo:PreservationOfEuclideaness}.
\item $f' \neq p$ are in general position in $\mathcal{O}^3 \cup p$  and are no coloops.
\item $(\mathcal{O}^3 \cup p,p,f')$ is Euclidean because of Corollary \ref{cor:lexExtStaysEucl2}, feasible if $\mathcal{O}^3$ is feasible and 
bounded because of Lemma \ref{lem:LexExtensionIsBounded}.
\item  $(\mathcal{O}^3 \cup p,p,f')$ has a topological sweep because of Theorem \ref{theo:topSweepExists}.
 \item Let $Z = (Z^0,\hdots, Z^s)$ be the cocircuits of the feasible region $P'(E_n)$ in $\mathcal{O}^3 \setminus f'$
 in sweep order.
\end{enumerate}

Let $\mathcal{F}'$ be  the face lattice of the feasible region of $\mathcal{O}^3 \setminus f'$. 
Because $ \rk(\mathcal{F}') = \rk(\mathcal{O}^3) = \rk-1$
the induction assumption is fullfilled and $Z$ induces a recursive atom ordering
of the cocircuits of the feasible region of $\mathcal{O}^3 \setminus f'$.

Lemma \ref{lem:EdgesAreCocircuits} yields that $\mathcal{F}'$ is isomorphic to $[Y^j,P(E_n)]$.
Let $f$ be the isomorphism defined in Lemma \ref{lem:EdgesAreCocircuits}.
Then $f(Z)$ is a recursive atom ordering
of the cocircuits of $[Y^j,P(E_n)]$.
In $Z$ the cocircuits with $f' = -$ come first. 
Then also the cocircuits with $f' = -$ come first in the recursive atom-ordering $f(Z)$  of $[Y^j,P(E_n)]$.
Because of Equation \ref{equ:atomsQj}
these are exactly the cocircuits of $Q_j$.
Be aware that the new order of the sweep can differ to the old one.
The definition of the recursive atom-ordering is completely fulfilled for $\mathcal{F}$ .
\end{proof}

We have the following Corollary.

\begin{corollary}\label{cor:VertexShellingExist}
 Let $\mathcal{O}$ be a Euclidean oriented matroid with $f \in \mathcal{O}$ lying in general position.
 Then each tope cell of $\mathcal{O} \setminus f$ is vertex-shellable.
\end{corollary}

\begin{proof}
With reorientation of some elements $\neq f$ of the matroid, the cocircuit graph $G_f$ stays the same. We can hence always assume that for the tope $X$ holds only $X_e \in \{0,+\}$ for all $e \in \mathcal{O} \setminus f$. Since $X$ is not empty it contains a cocircuit.
Let $\rk$ be the rank of $\mathcal{O}$. We pick independent elements $e_1, \hdots, e_r \in \mathcal{O} \setminus f$  and construct 
the lexicographical extension $\mathcal{O} \cup g = \mathcal{O}[e^+_1, \hdots, e^+_r]$. Because of Corollary \ref{cor:lexExtStaysEucl2} and Lemma \ref{lem:LexExtensionIsBounded}
the program $(\mathcal{O},g,f)$ is Euclidean, feasible and bounded, fullfilling the assumptions of Theorem \ref{theo:VertexShellingExist}.
\end{proof}

 \subsection{Vertex-shelling of the Edmonds-Mandel Face Lattice of a Euclidean Oriented Matroid}
 
 We conclude that the whole Euclidean oriented matroid is vertex-shellable.
 
\begin{theorem}\label{theo:VertexShellingExistWhole}
 Let $\mathcal{O}$ be a Euclidean oriented matroid of rank $\rk$ with $f \in \mathcal{O}$  in general position and  not a coloop.
Let $S = [e_1, \hdots, e_{r}]$ be an ordered set of independent elements of $\mathcal{O}$.
We define an ordering of the cocircuits of $\mathcal{O} \setminus f$ as follows.
\begin{enumerate}
\item Start with the unique cocircuit $X$ such that $X_{e_1} = \hdots = X_{e_{r-1}} = 0$ and $X_f = -$. Let $i = \rk-1$.
\item Continue with the cocircuits of index $i$ with $e_i = +$ and   $f = -$ and order them in the sweep order 
yielded by the Euclidean program  \\
$(\mathcal{O} / \{e_1, \hdots, e_{i-1}\}, e_i, f)$.
Then continue with the cocircuits of index $i$ with  $e_i = -$ and   $f = -$ and order them in the sweep order yielded by the same program of its negative cocircuits.
 Let $i = i-1$ and repeat (ii) until $i=0$.
\item    Let $i = i+1$. Continue with the cocircuits of index $i$ with $e_i = +$ and $f = +$ and order them in the sweep order
 yielded by the Euclidean program  \\
 $(\mathcal{O} / \{e_1, \hdots, e_{i-1}\}, e_i, f)$. 
 Then continue with the cocircuits of index $i$ with $e_i = -$ and $f = +$. Let $i = i+1$ and repeat (iii) until $i = \rk$.
\item Stop with the cocircuit $-X$.
\end{enumerate}
That is a vertex-shelling of the oriented matroid  $\mathcal{O} \setminus f$.
\end{theorem} 

\begin{proof}
  We proceed by induction on the rank $r\ge 2$ of $\mathcal{O}$. If
  the rank is two, there are only two cocircuits of index $2$ and all
  others are of index $1$. In this case the Euclidean program sorts
  the cocircuits allong the circle of its realization which clearly
  yields a shelling.
  
Let the Theorem hold for all oriented matroids of rank $\rk-1$ and let
$\mathcal{O}$ be a Euclidean oriented matroid of rank $\rk>2$.  Let $S =
(e_1, e_2, \hdots, e_r)$ be an ordered set of independent elements of
$\mathcal{O}$. We use the algorithm and obtain an ordering of the
cocircuits $X^1, \dots , X^m, Y^1, \hdots Y^n, Z^1, \dots Z^o$ of
$\mathcal{O} \setminus f$ such that $X^1, \dots , X^m$ are the
cocircuits of $\mathcal{O}$ with $e_1 = 0$ and $f = -$ and $Y^1, \dots
, Y^n$ the cocircuits with $e_1 \neq 0$ and $Z^1, \hdots , Z^o$ are
the cocircuits with $e_1 = 0$ and $f = +$.  The inductive assumption
yields that $X^1, \hdots, X^m, Z^1, \dots, Z^o$ is a vertex shelling
of $(\mathcal{O} / e_1) \setminus f$.

Let now $W^k$ be a cocircuit with $e_1 \neq 0$ and $f \neq 0$ and let $W^i$ be a cocircuit coming before 
$W^k$ in the sweep order. Then let $\mathcal{O}' = (\mathcal{O} \setminus f \cup f')$ be the oriented matroid yielded by moving the target function through $W^k$ (using the localization of Lemma \ref{lem:parallelExtensionEdmondsMandel}) hence we have  $W^k_{f'} = 0$. We have $W^i_{f'} = -$ in the case that $W^i_{e_1} = 0$ because of the parallelity of $f$ and $f'$ wrt. $e_1$  and otherwise
 because of Proposition \ref{prop:sweepOrientation2}
and because $W^i$ comes before $W^k$ in the sweep order.
We may assume $sep(W^i, W^k) = \emptyset$ otherwise  \cite{2}, Proposition 3.7.2 yields a cocircuit $W^j$ with $sep(W^j, W^k) = \emptyset$, $W^k \circ W^j \preceq W^k \circ W^i$ and $W^j_{f'} = -$.
We can go on with $W^j$ then because we want to prove that Theorem  \ref{theo:firstMainResult} (iii) holds.

Now let $\mathcal{O}'_T$ be be the oriented matroid $\mathcal{O}'$ reoriented
such that all elements of $supp(W^i \circ W^k) $ except $f'$ are positive. We still have $W^i_{f'} = -$ there.
Let $\mathcal{O}'' = (\mathcal{O}'_T  \cup p) = \mathcal{O}'_T[I^+]$ be the lexicographic extension with $I = (e_1, \hdots, e_r)$.
We obtain $W^i_p = W^k_p = +$ in $\mathcal{O}''$. 
Theorem \ref{:theo:PreservationOfEuclideaness} yields that the program $(\mathcal{O}'', p, f')$ is Euclidean with a topological sweep.
It is also bounded because of Lemma \ref{lem:LexExtensionIsBounded}.
The cocircuits $W^i$ and $W^k$ lie in the feasible region $(\mathcal{O}'', p, f')$.
Hence Theorem \ref{theo:firstMainResult}  (iii) yields a cocircuit $Z$ with $Z \circ Y^k$ being a conformal edge lying in $W^i \circ W^k$ and $Z_p = +$
directed from $Z$ to $W^k$ in $(\mathcal{O}'', p, f')$. Hence we have $Z_{f'} = -$.
Because $Z_p = +$ it is an old cocircuit (see Proposition \ref{prop:uniqueEdges2} ) existing already in $\mathcal{O}$.
If it holds $Z_{e_1} = 0$ we have $Z_f = -$ in $\mathcal{O}$ because of the parallelity of $f$ and $f'$ wrt. $e_1$.
Hence $Z$ comes before $W^k$ in the sweep order. 
If it holds $Z_{e_1} \neq 0$ we obtain $Dir(Z,W^k) = +$ in $(\mathcal{O}'_T, e_1, f')$ because of Proposition \ref{prop:basicDirectionProps} (iii)
hence also in $(\mathcal{O}, e_1, f)$ because of Proposition \ref{prop:parallelExtensionsDoNotChangeDirection}, Proposition \ref{prop:basicDirectionProps} (ii) and because reorienting elements from $E \setminus \{e_1,f\}$ does not
change the direction of any edge. Also in that case $Z$ comes before $W^k$ in the original sweep order.
Hence Theorem \ref{theo:firstMainResult} (iii) holds in $(\mathcal{O}, e_1, f)$.
We conclude that the formula in Theorem \ref{theo:firstMainResult} (iv) (see the proof there) holds.

Let now $W^k$ be a cocircuit with $W^k_{e_1} = 0$ and $W^k_f = +$.
Let $W^i$ be a cocircuit with $W^i_f = +$ and $W^i_{e_1} \neq 0$.
Then like before because of \cite{2}, Proposition 3.7.2 we may assume $sep(W^i,W^k) = \emptyset$ and the same Proposition yields a cocircuit $Z$ conformal to $W^k$ with $Z_{e_1} \neq 0$ and $Z \circ W^k$ being an edge
lying in $W^i \circ W^k$. It is clear that $Z$ is coming before $W^k$ in the sweep order and the formula in Theorem \ref{theo:firstMainResult} (iv) holds
again. 

Now we have only to consider the case that $W^k_{e_1} = 0$ and $W^i_{e_1} = 0$. 
Then for the cocircuit $Z$ must also hold $Z_{e_1} = 0$ we can look at the contraction $(\mathcal{O} / e_1, e_2,f)$.
The formula in Theorem \ref{theo:firstMainResult} (iv) holds again because of the induction assumption hence the non-recursive part (i) of Definition \ref{def:recAtomOrdering} 
of a vertex-shelling holds for $\mathcal{O} \setminus f$.

Now we have to show the recursive part of the definition. Let $Y^k$ be a cocircuit in $\mathcal{O} \setminus f$.
Let  $\mathcal{O}' = (\mathcal{O}  \cup p) = \mathcal{O}[I^+]$ be the lexicographic extension with $I = (e_1, \hdots, e_r)$. We have  $Y^k_p \neq 0$.
The program $(\mathcal{O}',p,f) $ is Euclidean. The ordering of the cocircuits given by our algorithm is a linear extension of the partial ordering induced by the graph $G_f$ of that program.
This is clear for the cocircuits having the same index because of Lemma \ref{lem:case3lexext} and with different indices because of Lemma \ref{lem:case6lexExt}.
Hence Theorem \ref{theo:topSweepExists} yields that we can move the target function  on $Y^k$ obtaining  $\mathcal{O}'' = \mathcal{O}' \setminus f \cup f'$ respecting our given ordering of the cocircuits.  
Then we go to $\mathcal{O}''' = \mathcal{O}'' \setminus supp(Y^k)$
which is of rank $\rk-1$, Euclidean because of Theorem
\ref{:theo:PreservationOfEuclideaness} and $\mathcal{O}''' \setminus
f'$ is isomorphic to $[Y^k, 1]$ because of Lemma
\ref{lem:EdgesAreCocircuits}. By induction assumption there is a
vertex shelling of $\mathcal{O}''' \setminus f'$ where the cocircuits
with $f' = -$ come first, but these are the cocircuits $W \setminus
supp(Y_k)$ with $W \circ Y^k$ being a conformal edge in $\mathcal{O}
\setminus f$ and $W$ is coming before $Y^k$ in our shelling order
because $W_{f'} = -$.
\end{proof}

 \section{concluding remarks}
 
 If we look carefully at the proofs of Theorem \ref{theo:VertexShellingExist} and Theorem \ref{theo:VertexShellingExistWhole}, we see that they remain valid
 if we assume the oriented matroid $\mathcal{O}$ only to be {\em Mandel}.
 Then we extend it to $\mathcal{O}' =  \mathcal{O} \cup f$ such that all programs $(\mathcal{O}',g,f)$ are Euclidean for all
 elements $g$ of the groundset of $\mathcal{O}$. Our last theorem says that $\mathcal{O} =  \mathcal{O}' \setminus f$ has a vertex-shelling then. Hence, an oriented matroid that is {\em Mandel} also has a vertex-shelling. 
 We are not aware, though, of an explicit example that is Mandel but not Euclidean.
We can immediately derive from our Theorem \ref{Theo:SecondMaintheorem2} that every Euclidean oriented matroid is {\em Mandel} but in general
not every oriented matroid needs to be {\em Mandel} (which was a conjecture made by A. Mandel in \cite{5}). J\"urgen Richter-Gebert found an
oriented matroid in \cite{Richter-Gebert} with a mutation-free element. That oriented matroid is a counterexample for being {\em Mandel} which was proven by K. Knauer
in \cite{Knauer}. It might be a candidate for an oriented matroid not being vertex-shellable at all. 
It requires further examinations to check that.

For subdivisions of zonotopes (which correspond to single-element liftings of realizable oriented matroids), 
Athanasiadis proved in \cite{Athanasiadis} that they are shellable (which means that the corresponding oriented matroid has a vertex shelling)
if the (realizable) oriented matroid corresponding to the zonotope is strongly Euclidean. 
Our last Theorem is a generalization of his result and 
it should be further examined what our results mean in the special case of zonotopal tilings.

%

\bibliographystyle{plain} 

\end{document}